\documentclass[sn-mathphys-num]{sn-jnl}

\usepackage{float}
\usepackage{graphicx}%
\usepackage{multirow}%
\usepackage{amsmath,amssymb,amsfonts}%
\usepackage{amsthm}%
\usepackage{mathrsfs}%
\usepackage{tabularx}
\usepackage[title]{appendix}%
\usepackage{xcolor}%
\usepackage{textcomp}%
\usepackage{manyfoot}%
\usepackage{booktabs}%
\usepackage{algorithm}%
\usepackage{algpseudocode}%
\usepackage{listings}%
\usepackage{subfig}

\theoremstyle{thmstyleone}%
\newtheorem{example}{Example}%
\theoremstyle{thmstyletwo}%
\theoremstyle{thmstylethree}%

\begin{document}

\title[Adaptive low-rank exponential integrators for DRE]{Adaptive low-rank exponential integrators for large-scale differential Riccati equation}


\author[1,2]{\fnm{Jingyi} \sur{Li}}\email{1249132@mail.dhu.edu.cn}

\author*[1]{\fnm{Dongping} \sur{Li}}\email{lidp@ccsfu.edu.cn}

\author[1]{\fnm{Hua} \sur{Yang}}\email{yanghzxcv@163.com}

\affil[1]{\orgdiv{School of Mathematics}, \orgname{Changchun Normal University}, \orgaddress{\street{677 Changji North Road}, \city{Changchun}, \postcode{130032}, \state{Jilin}, \country{People's Republic of China}}}

\affil[2]{\orgdiv{School of Mathematics and Statistics}, \orgname{Donghua University}, \orgaddress{\street{2999 Renmin North Road}, \city{Shanghai}, \postcode{201620}, \state{State}, \country{People's Republic of China}}}



\abstract{Matrix differential Riccati equation (DRE) typically exhibits transient and steady-state
phases, posing challenges for fixed-step time integration methods, which may lack accuracy during
transients or oversample in steady regimes. In this work, we propose adaptive low-rank matrix-valued exponential integrators for large-scale
stiff DRE. The methods combine embedded exponential Rosenbrock-type schemes and adaptive step-size control, enabling an automatic adjustment to the evolving
solution dynamics. This improves the accuracy during rapid transient phases while maintaining high accuracy in the steady state. Numerical experiments on benchmark problems demonstrate that the proposed adaptive integrators
consistently improve accuracy and computational efficiency compared with fixed-step low-rank
schemes.}




\maketitle

\section{Introduction}\label{sec:1}

We consider the matrix differential Riccati equation (DRE)
\begin{equation}\label{1.1}
\left\{
\begin{array}{l}
X'(t)=AX(t)+X(t)A^T+Q-X(t)GX(t)\equiv F(X),\\
X(0)=X_0,
\end{array}
\right.
\end{equation}
where $A\in\mathbb{R}^{N\times N}$, and the matrices $Q$, $G$, and $X_0$ are symmetric. Under mild assumptions,  the exact solution of \eqref{1.1} converges to its equilibrium $X^*$ as $t\longrightarrow \infty $, which satisfies the algebraic Riccati equation (ARE)
\begin{equation}\label{1.2}
AX^*+X^*A^T+Q-X^*GX^*=0.
\end{equation}

The DRE arises in a broad range of applications, including optimal control, filtering, model reduction, and differential games \cite{Abou03, Antoulas2005, Benner17, Reid1972}. In practice, such problems are often high-dimensional, e.g., from spatial discretizations of PDEs, resulting in large, sparse system matrix $A$. Efficient time integration is crucial both in terms of computation and memory. Moreover, the coexistence of fast and slow dynamical modes typically renders DRE stiff, posing additional challenges for  efficient and accurate numerical methods.

Numerous numerical approaches have been developed for solving DRE; see, e.g., \cite{choi1990, dieci92, dieci94, mena07, benner13}. Classical methods, however, are generally limited to small- or medium-scale problems and become computationally prohibitive for large-scale systems. In recent years, considerable attention has been devoted to large-scale DRE, motivated by the observation that in many applications the solution matrix exhibits rapidly decaying singular values; see \cite{stillfjord18B} for theoretical analysis.  Building on this idea, several low-rank integrators have been proposed for large-scale stiff DRE. These include matrix versions of BDF and Rosenbrock methods \cite{lang15,Benner18A}, Peer methods \cite{Benner18B}, splitting methods \cite{Stillfjord15, Stillfjord18a, Mena2018, Ostermann2019}, and Krylov subspace methods \cite{Koskela2020,Simoncini19}. More recently, low-rank exponential integrators have been proposed and shown to be particularly effective for large-scale stiff problems \cite{Li2021, Li2025}.

Despite these advances, most existing integrators are built upon fixed time-stepping strategies, which may limit either accuracy or efficiency. The solution of a DRE typically exhibits pronounced multiscale behavior: a short transient phase dominated by fast dynamics, followed by a slowly evolving steady-state regime. A fixed time step is inherently incapable of resolving these qualitatively different dynamical features in a balanced manner.
For instance, splitting-based methods often capture the initial dynamics effectively. However, their accuracy may deteriorate over long-time integration. Exponential integrators, on the other hand, are well known for their ability to maintain high accuracy near steady state. Yet, when employed with fixed time steps, they may fail to efficiently resolve rapid transient behavior unless prohibitively small steps are chosen. Consequently, fixed-step strategies impose an intrinsic trade-off between accuracy and efficiency: step sizes small enough to resolve transient dynamics lead to unnecessary computational effort in the steady-state regime, whereas larger steps reduce cost at the expense of accuracy. This limitation is not incidental but structural, and it restricts the practical performance of existing low-rank integrators for large-scale stiff DRE.

These considerations strongly motivate the development of adaptive time integration techniques. In this work, we propose adaptive low-rank matrix-valued exponential integrators for large-scale stiff DRE. By combining embedded exponential Rosenbrock-type schemes with efficient low-rank implementations and adaptive step-size control, the proposed methods automatically adjust to the evolving solution dynamics. This leads to improved accuracy and computational efficiency across both transient and steady-state regimes. Numerical experiments on benchmark problems demonstrate that the proposed adaptive integrators provide a competitive and scalable alternative to existing fixed-step low-rank methods.

The remainder of the paper is organized as follows.
Section~\ref{sec:2} reviews exponential Rosenbrock-type methods and discusses their low-rank implementation for large-scale DRE.
Section~\ref{sec:3} introduces the proposed adaptive low-rank exponential Rosenbrock integrators, focusing on error estimation and practical step-size control within a low-rank framework.
Numerical experiments in Section~\ref{sec:4} illustrate the accuracy and efficiency of the proposed methods on benchmark problems.
Finally, Section~\ref{sec:5} concludes the paper.

\section{Constant-step low-rank exponential Rosenbrock integrators for DRE}\label{sec:2}

Exponential integrators have a long history and have been shown to be competitive with classical implicit methods for the time integration of semi-linear stiff problems. A key feature of this class of methods is that the dominant linear part of the system is treated exactly, while the nonlinear part is handled explicitly. This strategy leads to favorable stability properties and enables explicit time integration even for stiff problems. Over the past decades, a large variety of exponential integrators have been developed, analyzed, and implemented; see, e.g., \cite{hochbruck97,kassam05,MH2006,Tokman2006,MH2011,Tokman11,Luan2013,loffeld13,Luan2016}. For a comprehensive overview, we refer to the survey \cite{Hochbruck10} and the references therein.

In this work, we focus on exponential Rosenbrock-type methods, originally introduced in \cite{hochbruck09}. These methods form an important subclass of exponential integrators based on a continuous linearization of the vector field along the numerical trajectory. Compared with general exponential integrators, exponential Rosenbrock methods offer improved stability behavior and significantly simplified order conditions. They are particularly attractive for the DRE, since the quadratic structure of the right-hand side allows the linearization of the vector field to be computed analytically. As a result, the required Jacobian-related operations can be evaluated efficiently without additional approximation or preprocessing. This property makes exponential Rosenbrock methods well suited for large-scale problems.

\subsection{Matrix-valued exponential Rosenbrock integrators for DRE}

We next review matrix-valued exponential Rosenbrock integrators for the DRE \eqref{1.1}, following the low-rank formulation in \cite{Li2021}. This review is included for completeness and to establish the notation and algorithmic components required for the development of adaptive low-rank schemes in the next section.

Let $X_n \approx X(t_n)$ denote the numerical solution at time $t_n$. Exploiting the structure of the DRE, the right-hand side in \eqref{1.1} can be decomposed into a linear and a nonlinear part,
\begin{equation}\label{2.2}
X'(t)=F(X(t))=\mathcal{L}_n[X]+\mathcal{N}_n(X),     
\end{equation}
where $\mathcal{L}_n(X)$ denotes the Fr\'{e}chet derivative of $F$ evaluated at $X_n$, and $\mathcal{N}_n(X)$ is the nonlinear remainder. Specifically,
\begin{eqnarray}\label{2.3}
\mathcal{L}_n[X]=A_nX+XA_n^T,\qquad \mathcal{N}_n(X)=F(X)-\mathcal{L}_n[X]
\end{eqnarray}
with $A_n=A-X_nG$.

 By the variation-of-constants formula, the exact solution of \eqref{2.2} satisfies
\begin{equation}\label{2.4}
X(t_{n}+h_n)=e^{h_n\mathcal{L}_n}\left[X(t_n)\right]+h_n\int_{0}^{1}e^{(1-s)h_n\mathcal{L}_n}\left[\mathcal{N}_n(X(t_n+sh_n))\right]ds.
\end{equation}
Exponential Rosenbrock methods are obtained by approximating the nonlinear term
$\mathcal{N}_n(X(t_n+s h_n))$ in \eqref{2.4} by suitable interpolation polynomials based on intermediate stage values. This leads to the following general $s$-stage exponential Rosenbrock integrators:
\begin{equation}\label{2.5}
\begin{cases}
X_{ni}=e^{c_ih_n\mathcal{L}_n}[X_n]+h_n\sum\limits^{i-1}_{j=1}a_{ij}(h_n\mathcal{L}_n)[\mathcal{N}_n(X_{nj})],~~1\leq i\leq s,\\[2ex]
X_{n+1}=e^{h_n\mathcal{L}_n}[X_n]+h_n\sum\limits^{s}_{i=1}b_{i}(h_n\mathcal{L}_n)[\mathcal{N}_n(X_{ni})].
\end{cases}
\end{equation}
Here, $h_n>0$ denotes the time step size at $t_n$, and $c_i$ are the stage 
nodes, and $X_{ni}$, $X_{n+1}$ represent the numerical approximations at times $t_n+c_i h_n$ and $t_{n+1}$, respectively. The coefficient functions $a_{ij}(z)$ and $b_i(z)$ are linear combinations of the $\varphi$-functions evaluated at $c_i z$ and $z$, respectively. These functions are defined as
\begin{equation}
\varphi_0(z) = e^z, \quad \varphi_k(z) = \int_0^1 e^{(1-\theta)z} \frac{\theta^{k-1}}{(k-1)!} \,d\theta, \quad k \geq 1.
\end{equation}

Imposing 
\begin{equation}
\sum_{i=1}^{s} b_i(z)=\varphi_1(z), \qquad
\sum_{j=1}^{i-1} a_{ij}(z)=c_i\,\varphi_1(c_i z),
\quad 1\le i\le s,
\end{equation}
ensures that \eqref{2.5} preserves the equilibrium solution of the DRE. The scheme can then be reformulated in the form
\begin{equation}\label{2.6}
\begin{cases}
X_{ni}&=X_n+c_ih_n \varphi_1\left(c_ih_n\mathcal{L}_n\right)[F(X_n)]+h_n\sum\limits^{i-1}_{j=2}a_{ij}(h_n\mathcal{L}_n)[\mathcal{D}_{nj}],~~2\leq i\leq s,\\
X_{n+1}&=X_n+h_n \varphi_1\left(h_n\mathcal{L}_n\right)[F(X_n)]+h_n\sum\limits^{s}_{j=2}b_{j}(h_n\mathcal{L}_n)[\mathcal{D}_{nj}],
\end{cases}
\end{equation}
where 
\begin{equation}\label{2.6c}
\mathcal{D}_{nj}=\mathcal{N}_n(X_{nj})-\mathcal{N}_n(X_{n}), \quad 2\leq j\leq s.
\end{equation}
 The simplest scheme is the exponential Rosenbrock–Euler scheme,
\begin{equation}\label{2.7}
X_{n+1}=X_n+h_n\varphi_1(h_n\mathcal{L}_n)[F(X_n)],
\end{equation}
which is a second-order accurate method requiring only a single evaluation of an operator-valued $\varphi$-function at each time step.
The construction and convergence analysis of matrix-valued exponential Rosenbrock integrators closely parallel those of the corresponding vector-valued methods and are therefore not repeated here; see \cite{hochbruck09,Luan2013,Luan2014} for details.

It is important to note that, in contrast to vector-valued exponential integrators, the matrix-valued setting involves $\varphi$-functions of Lyapunov operators. The efficient evaluation of these operator functions is therefore crucial for the overall performance of the matrix-valued integrators. Both full-rank and low-rank techniques for this purpose have been investigated in \cite{Li2023,Li2025}.

\subsection{Low-rank implementation of matrix-valued exponential Rosenbrock integrators}
In many large-scale applications, the DRE \eqref{1.1} exhibits an intrinsic low-rank structure, and its solution typically remains numerically low-rank over time. Explicitly forming full matrix approximations would therefore lead to dense computations and prohibitive storage costs. To reduce the computational cost, scheme \eqref{2.6} must be implemented directly in terms of low-rank level. Specifically, the numerical solution is represented in factored form, and all matrix operations are performed in a low-dimensional subspace. This allows the overall computational complexity to scale linearly with the state dimension. 

In the following, we briefly outline the low-rank formulation and highlight the key computational building blocks that will be used in the adaptive integrators developed in Section~\ref{sec:3}.

From now on, we assume that the coefficient matrices $Q$, $G$, and the initial value $X_0$ admit low-rank factorizations
\begin{equation}\label{2b.1}
\begin{aligned}
Q&=C^TC, && C\in\mathbb{R}^{p\times N},\vspace{1.5ex}\\
G&=BB^T,&& B\in\mathbb{R}^{N\times q},\vspace{1.5ex}\\
X_0&=L_0D_0L_0^T,&& L_0\in\mathbb{R}^{N\times r},\; D_0\in\mathbb{R}^{r\times r},
\end{aligned}
\end{equation}
where $p,q,r\ll N$.

Assume that the numerical solution at time $t_n$ admits the factorization
\[
X_n=L_nD_nL_n^T,\qquad L_n\in\mathbb{R}^{N\times r_n},\; D_n\in\mathbb{R}^{r_n\times r_n}.
\]
with $r_n\ll N$.
Our objective is to update the low-rank factors $L_{n+1}$ and $D_{n+1}$ of $X_{n+1}$ directly, without explicitly forming the full matrix  $X_{n+1}$, thereby reducing the computational cost.

Using these representations, the right-hand side $F(X_n)$ in \eqref{1.1} admits the factorization
\begin{equation}\label{2b.2}
F(X_n)=\widetilde L_n \widetilde D_n \widetilde L_n^T,
\end{equation}
with
\[
\widetilde L_n=[\,C^T,\; AL_n,\; L_n\,],
\]
and
\[
\widetilde D_n=
\begin{pmatrix}
I_p & 0 & 0 \\
0 & 0 & D_n \\
0 & D_n & -(D_nL_n^TB)(D_nL_n^TB)^T
\end{pmatrix}.
\]
The matrix $\widetilde L_n$ typically contains substantially more columns than $L_n$, and its column dimension may exceed the numerical rank of $F(X_n)$. Consequently, without additional treatment, the rank of the intermediate factorizations may grow unnecessarily, leading to increased computational and memory costs.
To reduce the computational cost and control rank growth, a compression procedure -- such as truncated SVD or QR-based truncation \cite{lang15}—is applied whenever new low-rank factors are generated. This procedure eliminates redundant directions in the factorization, keeps the effective rank bounded during the time integration, and thereby preserves the efficiency and saves memory.

For the term
\begin{equation}\label{2b.3}
\begin{aligned}
\mathcal{D}_{nj}
&=\mathcal{N}_n(X_{nj})-\mathcal{N}_n(X_n)\\
&=X_nGX_{nj}+X_{nj}GX_n-X_{nj}GX_{nj}-X_nGX_n,
\end{aligned}
\end{equation}
 it can be expressed compactly as
\begin{equation}\label{2b.4}
\mathcal D_{nj}=-K_{nj} G K_{nj}, \quad K_{nj}=X_{nj}-X_n.
\end{equation}
If the stage value $X_{nj}$ admits a low-rank factorization
\begin{equation}\label{2b.5}
 X_{nj}=L_{nj} D_{nj} L_{nj}^T,
\end{equation}
then the difference $ K_{nj}$  can be written as $K_{nj}=U_{nj}T_{nj}U_{nj}^T$  with $U_{nj}=[L_n, L_{nj}]$ and $T_{nj}=\operatorname{blkdiag}(-D_n, D_{nj})$.
Consequently, the term $\mathcal D_{nj}$ itself admits a low-rank factorization
\begin{equation}\label{2b.6}
\mathcal D_{nj}=\widehat{L}_{nj}\widehat{\Gamma}_{nj}\widehat{L}_{nj}^T,
\end{equation}
with
\begin{equation}\label{2b.7}
\widehat{L}_{nj}=U_{nj} \quad \text{and} \quad \widehat{\Gamma}_{nj}=-(T_{nj}U_{nj}^TB)(T_{nj}U_{nj}^TB)^T.
\end{equation}

Since the quantities 
$a_{ij}(h_n\mathcal{L}_n)[\mathcal{D}_{nj}]$ and 
$b_j(h_n\mathcal{L}_n)[\mathcal{D}_{nj}]$ in \eqref{2.6} 
can be expressed as linear combinations of 
$\varphi_j(c \mathcal{L}_n)[\mathcal{D}_{nj}]$, 
with $c=c_ih_n$ or $c=h_n$, 
the stage values $X_{ni}$ and the updated solution $X_{n+1}$ 
can also be written as linear combinations of these operator-valued 
$\varphi$-functions. Consequently, the computational bottleneck of the exponential Rosenbrock scheme \eqref{2.6} 
reduces to the efficient evaluation of 
$\varphi_j(c h_n\mathcal{L}_n)[\mathcal{D}_{nj}]$.

In this work, the low-rank factorizations of
$\varphi_j(c h_n\mathcal{L}_n)[\mathcal{D}_{nj}]$
are computed using the recursive low-rank solver \texttt{recurLrlyap} \cite{Li2025}. This approach exploits the structure of the low-rank structure of $\mathcal{D}_{nj}$ and generates the factors in a recursive manner. Its computational complexity scales with the rank of $\mathcal{D}_{nj}$ and the sparsity of the involved operators. As a result, the overall complexity of the proposed exponential integrators is governed by low-rank linear algebra operations, making the method particularly suitable for large-scale problems.

Once low-rank factorizations of these operator $\varphi$-functions are available, the remaining of the scheme \eqref{2.6}  can be carried out
entirely in the low-rank format. In particular,
both the stage values $X_{ni}$ and the updated solution $X_{n+1}$ 
naturally admit low-rank representations that can be constructed directly from the corresponding linear combinations.
This avoids the formation of dense matrices and ensures that the computational and memory costs depend primarily on the ranks of the intermediate quantities.

To make the presentation more explicit, consider a generic representation of the updated solution.
Suppose that $X_{n+1}$ can be written in the form
\begin{equation}\label{2b.8}
X_{n+1}=\sum_{k=1}^{m} g_{n,k} \widetilde{L}_{n,k} \widetilde{D}_{n,k}\widetilde{ L}_{n,k}^T.
\end{equation}
Then a low-rank $LDL^T$ factorization of $X_{n+1}$ can be obtained by concatenating the basis matrices,
\[
L_{n+1} = [\widetilde{L}_{n,1}, \widetilde{L}_{n,2}, \ldots, \widetilde{L}_{n,m}], 
\qquad 
D_{n+1} = \operatorname{blkdiag}(g_{n,1} \widetilde{D}_{n,1}, g_{n,2} \widetilde{D}_{n,2}, \ldots, g_{n,m} \widetilde{D}_{n,m}),
\]
which provides the factorization of $X_{n+1}$ without explicitly forming full matrix.

\section{Adaptive low-rank exponential Rosenbrock integrators for DRE}\label{sec:3}
In this section, building on the low-rank formulation introduced in Section \ref{sec:2}, we present adaptive exponential integration framework tailored for large-scale DRE. In particular, we develop embedded error estimation and adaptive time-stepping strategies within low-dimensional subspaces, so that the low-rank structure is preserved throughout the integration.

\subsection{Embedded exponential Rosenbrock integrators pair}
Let $X_{n+1}$ denote the numerical solution produced by the exponential Rosenbrock scheme \eqref{2.6} of order $q$. 
To enable adaptive step-size control, we consider an embedded approximation
\begin{equation}\label{3.1}
\overline{X}_{n+1}=X_n+h_n \varphi_1\left(h_n\mathcal{L}_n\right)[F(X_n)]+h_n\sum\limits^{s}_{j=2}\overline{b}_{j}(h_n\mathcal{L}_n)[\mathcal{D}_{nj}],
\end{equation}
which sharing the same internal stages $X_{ni}$ but employs different weights $\overline{b}_{j}$ in output stage. 
Assume the embedded method has order $p\leq q$.  The corresponding local error estimate is given by
\begin{equation}\label{3.2}
E_{n+1}=h_n\sum\limits^{s}_{i=2} \left(b_{i}(h_n\mathcal{L}_n)-\overline{b}_{i}(h_n\mathcal{L}_n)\right) \left[\mathcal{D}_{ni}\right],
\end{equation}
and satisfies the asymptotic relation
\begin{equation*}
\|E_{n+1}\| = c_{n+1} h_n^{p+1} + \mathcal{O}(h_n^{p+2}).
\end{equation*}

It is important to emphasize that, within the proposed low-rank framework, the local error estimator admits a low-rank factorization
\begin{equation*}
 E_{n+1} = L_{E_n}D_{E_n}{L_{E_n}}^T.
\end{equation*}
The Frobenius norm of $E_{n+1}$ required for step-size control can be computed via
\begin{equation*}
 \|E_{n+1}\| = \text{trace}\left(({L_{E_n}}^T L_{E_n}D_{E_n})^2\right)^{1/2}.
\end{equation*}
In particular, if $L_{E_n}$ has orthogonal columns, i.e., ${L_{E_n}}^T L_{E_n}=I$, then
\begin{equation*}
 \|E_{n+1}\| = \|D_{E_n}\|.
\end{equation*}
Hence, the Frobenius norm of the local error can be computed directly from the small matrix $D_{E_n}$,
 without assembling the full matrix $E_{n+1}$.
This significantly reduces the computational cost of the error estimation step.

For brevity, we write $b_i = b_i(z)$, $\overline{b}_{i}=\overline{b}_{i}(z)$, and $a_{ij} = a_{ij}(z)$. 
The embedding scheme can be denoted by the reduced Butcher tableau:
\begin{equation*}
\begin{array}{c|cccc}
c_2 &  &\\
c_3 & a_{32} &\\
\vdots & \vdots & \ddots\\
c_s & a_{s2} & \ldots & a_{s,s-1} \\
\hline 
& b_2 & \ldots & b_{s-1} & b_s\\
& \overline{b}_2 & \ldots &\overline{ b}_{s-1} & \overline{b}_s\\
\end{array}
\end{equation*}

In the numerical experiments of Section~\ref{sec:4}, we adopt two representative embedded exponential Rosenbrock-type schemes for adaptive time stepping. The corresponding coefficients and detailed formulations are given below.

$\bullet$ \texttt{exprb32}. The scheme is of third order, which combines a third-order exponential Rosenbrock integrator with a second-order embedded approximation \cite{hochbruck09}. Its Butcher tableau is given by
\begin{equation*}
\begin{array}{c|cc}
c_2 &   \\
\hline & b_2\\
 & \overline{b}_2  \\
\end{array}~~~~=~~\begin{array}{c|cc}
1  & \\
\hline  & 2\varphi_{3}\\
 & 0
\end{array}
\end{equation*}
The method consists of one internal stage,
More explicitly, the scheme reads
\begin{equation}\label{2.8}
\begin{cases}
X_{n2}&=X_n+h_n\varphi_1(h_n\mathcal{L}_n)[F(X_n)],\\[2ex]
X_{n+1}&=X_n+h_n\varphi_1(h_n\mathcal{L}_n)[F(X_n)]+2h_n\varphi_3(h_n\mathcal{L}_n)\left[\mathcal{D}_{n2}\right],\\
\overline{X}_{n+1}&=X_{n2}.

\end{cases}
\end{equation}

$\bullet$ \texttt{exprb43}. The scheme is a fourth-order exponential Rosenbrock scheme with a third-order embedding \cite{Luan2016}. 
 Its coefficients are
\begin{equation*}
\begin{array}{c|ccc}
c_2 &  & \\
c_3 & a_{32} & \\
\hline & b_2 & b_3\\
& \overline{b}_2 & \overline{b}_3\\
\end{array}~~~~=~~\begin{array}{c|ccc}
\frac{1}{2} &  & \\
1  &  & \\
\hline  & 16\varphi_3-48\varphi_4 &-2\varphi_{3}+12\varphi_4\\
&16\varphi_3 &-2\varphi_3
\end{array}
\end{equation*}
 The method proceeds as follows
\begin{equation}
\begin{cases}
X_{n2}&=X_n+\frac{1}{2}h_n\varphi_1(\frac{1}{2}h\mathcal{L}_n)[F(X_{n})]\\
X_{n3}&=X_n+h\varphi_1\left(h\mathcal{L}_n\right)[F(X_{n})]\\
\overline{X}_{n+1}&=X_n+h\varphi_1(h\mathcal{L}_n)[F(X_{n})]+16h_n\varphi_3(h_n\mathcal{L}_n)[\mathcal{D}_{n2}] -2h_n\varphi_3(h_n\mathcal{L}_n)[\mathcal{D}_{n3}]\\  
X_{n+1}&=\overline{X}_{n+1}-48h_n\varphi_4(h_n\mathcal{L}_n)[\mathcal{D}_{n2}]+12h_n\varphi_4(h_n\mathcal{L}_n)[\mathcal{D}_{n3}].
\end{cases}
\end{equation}

For the embedded exponential Rosenbrock schemes \texttt{exprb32} and \texttt{exprb43}, the local error estimators admit explicit representations in terms of higher-order Lyapunov operator $\varphi$-functions. In particular, the error estimates are given by 
\begin{equation*}
2h_n\left\|\varphi_3(h_n\mathcal{L}_n)\left[\mathcal{D}_{n2}\right]\right\| \quad \text{and} \quad h_n\left\|\varphi_4(h_n\mathcal{L}_n)[-48\mathcal{D}_{n2}+12\mathcal{D}_{n3}]\right\|,
\end{equation*}
 respectively.
 As a result,  the local error estimator can be evaluated at the factor level using the recursive low-rank solver \texttt{recurLrlyap}.

\subsection{Step-size selection}
We employ a standard adaptive step-size  strategy following \cite[Chapter II.4]{hairer93}.
Let 
\begin{equation*}
\texttt{Tol}=\texttt{Atol}+\max\{\|X_n\|,\|X_{n+1}\|\}\cdot\texttt{Rtol},
\end{equation*}
 where \texttt{Atol} and \texttt{Rtol} denote the prescribed absolute and relative tolerances, respectively.
If $\|E_{n+1}\| \leq \texttt{Tol}$, the computed step is accepted and the next step size is chosen as
\begin{equation}
h_{n+1}= \min\left\{\delta_{max}, \sigma_1 \left( \frac{\texttt{Tol}}{\|E_{n+1}\|} \right)^{\frac{1}{p+1}}\right\}h_n. 
\end{equation}
Otherwise, the step is rejected and recomputed with
\begin{equation}
h_n= \max\left\{\delta_{\min}, \sigma_2 \left( \frac{\texttt{Tol}}{\|E_{n+1}\|} \right)^{\frac{1}{p+1}}\right\}h_n. 
\end{equation}
Here, $\sigma_1$ and $\sigma_2$ are safety factors, while  $\delta_{max}$ and $\delta_{min}$ restrict the maximal increase and minimal decrease of the step size  in order to avoid overly aggressive changes. This implies that the new step size  can be at most $\delta_{max}$ times and at least $\delta_{min}$ times the current step size.
In our implementation, we set $\sigma_1=0.9, \delta_{max}=1.5$,  $\sigma_2=0.5$, and $\delta_{min}=0.1$.

The selection of an appropriate initial step size typically relies on a reliable estimate of the local error. However, for general high-order exponential integrators, the explicit computation of the local truncation error is often complicated and computationally expensive. To determine a suitable initial step size within our adaptive framework, we instead exploit the local truncation error of the second-order exponential Rosenbrock–Euler scheme. This method provides a robust, inexpensive, and sufficiently accurate error indicator for the initialization phase. Moreover, within the proposed low-rank framework, all quantities required in this procedure can be evaluated directly in factored form.

The local truncation error (LTE) of the exponential Rosenbrock–Euler scheme satisfies the asymptotic relation 
\begin{equation}
\text{LTE}\approx Ch_0^3\|F(X_0)GF(X_0)\|,
\end{equation}
where $C$ is a problem-dependent constant. The quantity $\|F(X_0)GF(X_0)\|$ can be computed in a low-dimensional subspace by exploiting the low-rank structures of both $F(X_0)$ and $G=BB^T$. By equating this estimate with the tolerance
\begin{equation}
Ch_0^3\|F(X_0)GF(X_0))\|=\texttt{Tol},\quad \quad \texttt{Tol}=\texttt{Atol}+\|X_0\|\cdot\texttt{Rtol},
\end{equation}
we obtain the initial step size
\begin{equation}
h_0=\theta \left(\frac{\texttt{Tol}}{\|F(X_0)GF(X_0))\|}\right)^{1/3},
\end{equation}
where $\theta =C^{-1/3}$. 

In practice, $\theta$ is treated as a safety factor and set to  $\theta=0.1$.
In the implementation, the exponent $1/3$ can be replaced by $1/(p+1)$, where $p$ denotes the order of the time integrator.
Numerical experiments demonstrate that this initialization strategy is effective for stiff and large-scale DRE, where the solution often undergoes a highly dynamic transient phase, and reliable step-size selection is crucial for overall efficiency.

\section{Numerical experiments}\label{sec:4}

In this section, we report a set of numerical experiments designed to evaluate the accuracy and efficiency of the adaptive low-rank exponential Rosenbrock integrators, \texttt{exprb32} and \texttt{exprb43}. 
Their performance is compared with several established low-rank time integrators for large-scale DRE. 
In particular, we consider fixed time-step exponential Rosenbrock integrators of orders 2 and 3 (denoted by \texttt{exprb2} and \texttt{exprb3}) \cite{Li2021}, the second-order Rosenbrock method (denoted by \texttt{Ros2}) \cite{lang15}, and a fourth-order splitting scheme \cite{Stillfjord18a}. 
State-of-the-art implementations of the latter two methods are provided in the M.E.S.S. toolbox \cite{SaaKB21mmess}. 

All experiments were performed in MATLAB on a desktop computer equipped with an Intel Xeon G5128 processor (2.1\,GHz) and 64\,GB of RAM. 
Unless otherwise specified, the accuracy is measured by the relative error at the final time in the Frobenius norm.

\begin{figure}[H]
\begin{minipage}{0.5\linewidth}
\centering
\includegraphics[width=6.5cm,height=5cm]{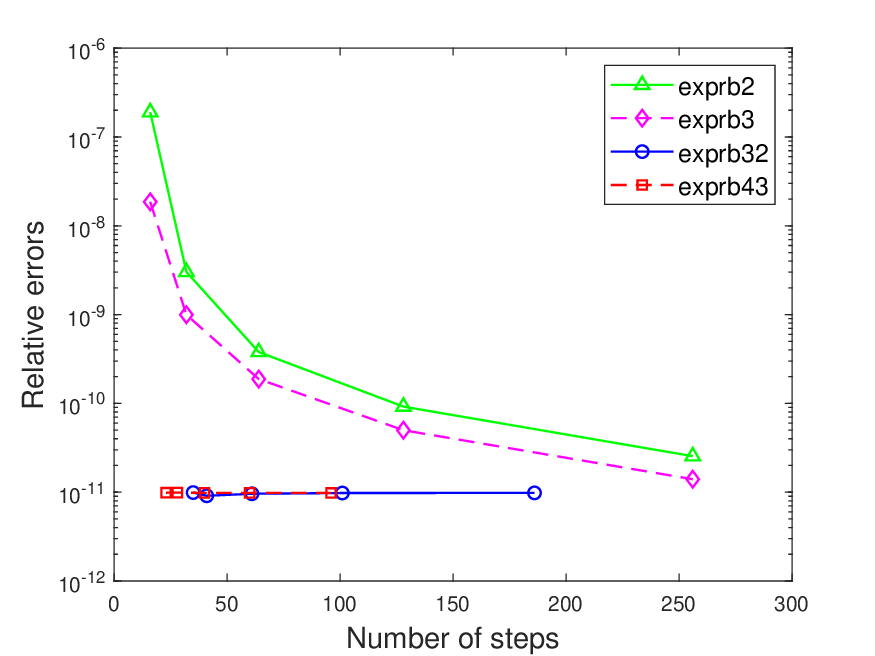}\\
\end{minipage}
\mbox{\hspace{-0.5cm}}
\begin{minipage}{0.5\linewidth}
\centering
\includegraphics[width=6.5cm,height=5cm]{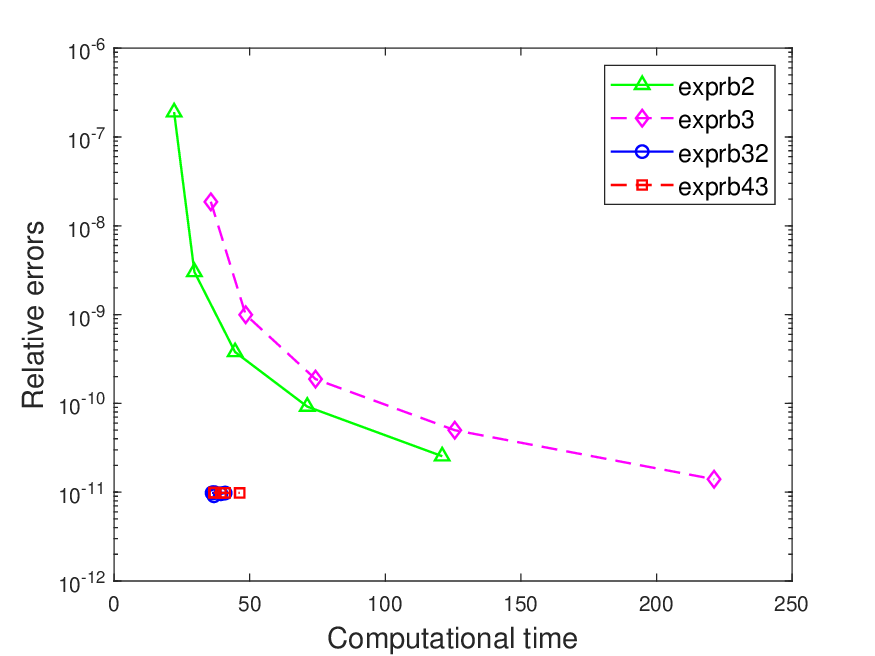}\\
\end{minipage}
\caption{Comparison between the adaptive and non-adaptive low-rank exponential integrators when integrating DRE over [0, 0.1] in Experiment \ref{exam1}. Left: Relative errors as a function of number of steps. Right: Relative errors as computational time.}\label{fig4.1}
\end{figure}

\begin{figure}[H]
\begin{minipage}{0.5\linewidth}
\centering
\includegraphics[width=6.5cm,height=5cm]{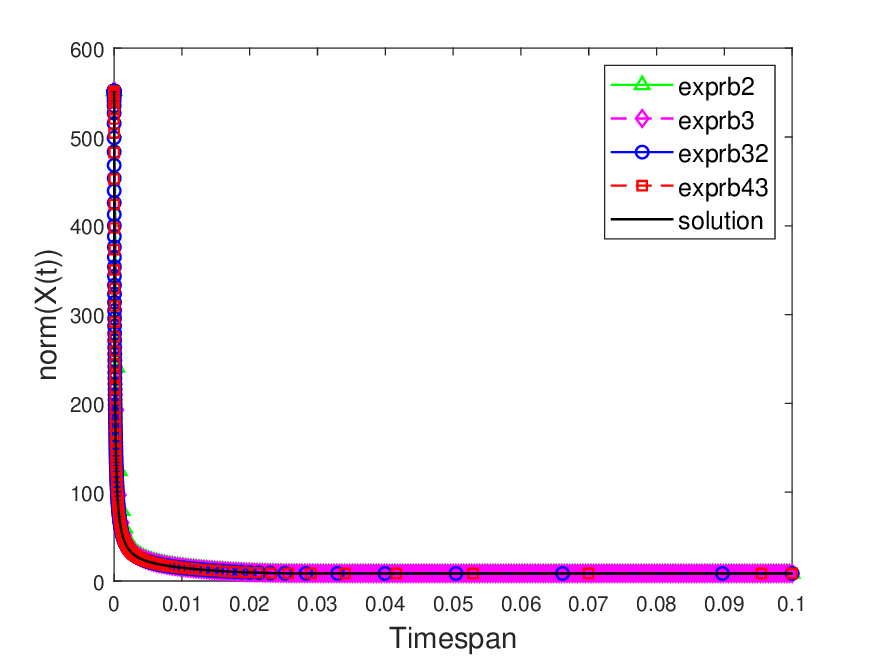}\\
\end{minipage}
\mbox{\hspace{-0.5cm}}
\begin{minipage}{0.5\linewidth}
\centering
\includegraphics[width=6.5cm,height=5cm]{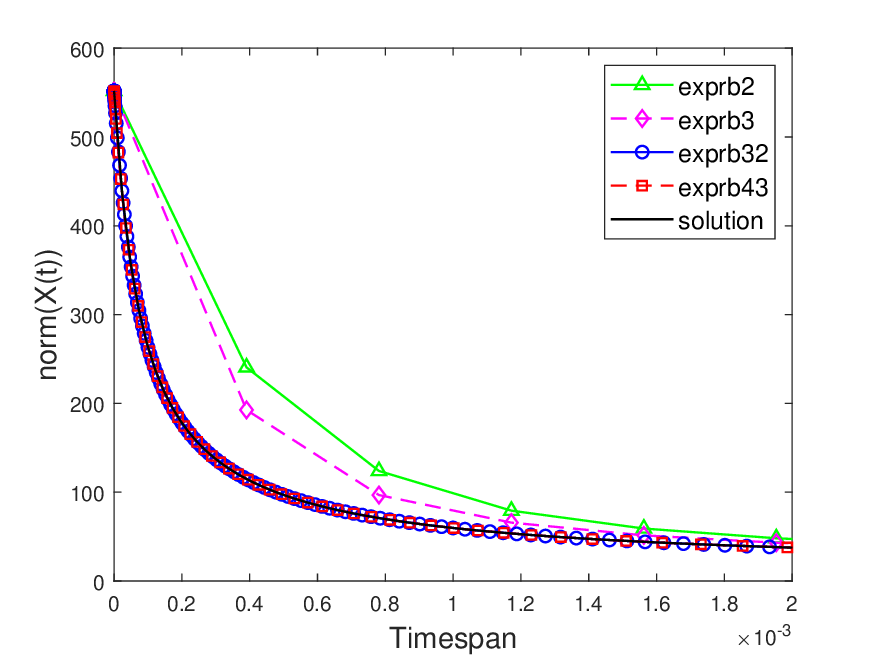}\\
\end{minipage}
\caption{Evolution of the reference, adaptive, and non-adaptive low-rank exponential integrators for DRE in Experiment \ref{exam1}. Left: Results over [0, 0.1]. Right: Zoom of the left panel.}\label{fig4.2}
\end{figure}

\begin{figure}[H]
\begin{minipage}{0.5\linewidth}
\centering
\includegraphics[width=6.5cm,height=5cm]{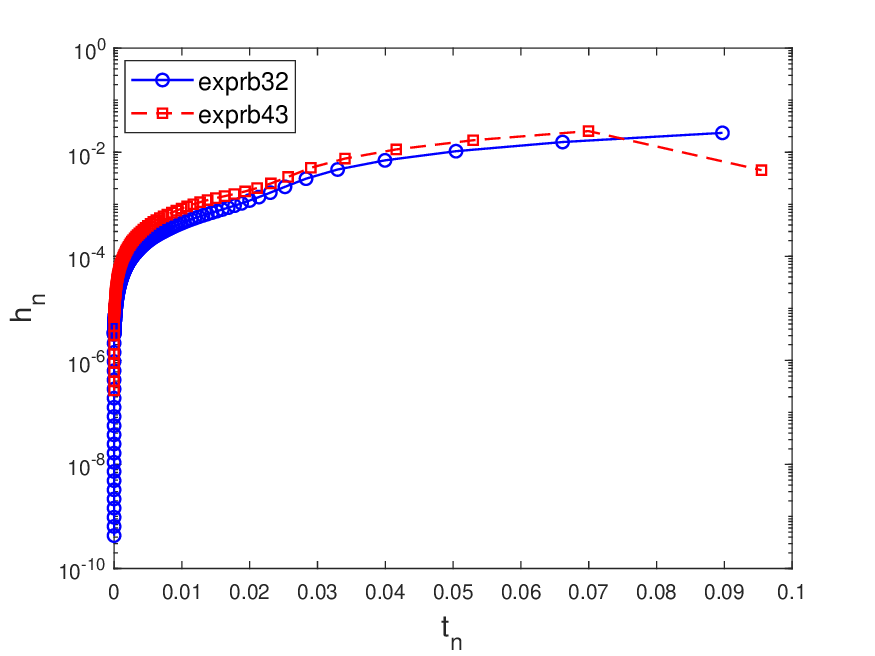}\\
\end{minipage}
\mbox{\hspace{-0.5cm}}
\begin{minipage}{0.5\linewidth}
\centering
\includegraphics[width=6.5cm,height=5cm]{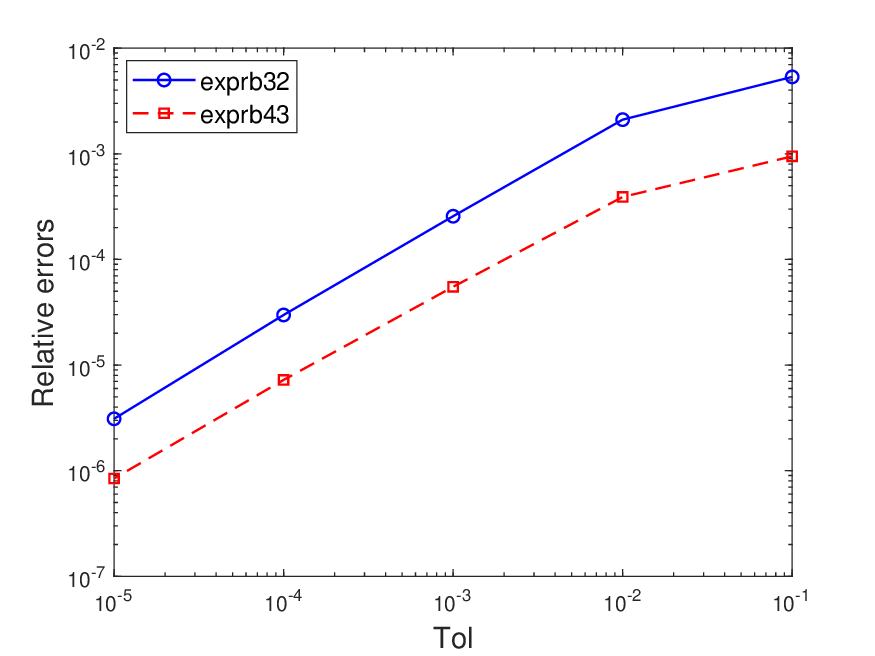}\\
\end{minipage}
\caption{Adaptive time step sizes and relative errors for varying tolerances for the DRE in Experiment \ref{exam1}. Left: Semi-log plot of adaptive time step sizes. Right: Relative errors at $t = 0.002$ as a function of prescribed tolerance.} \label{fig4.3}
\end{figure}

\begin{example}\label{exam1}
In this example, we compare the performance of the proposed adaptive and non-adaptive low-rank exponential integrators. The test problem is a widely used benchmark DRE taken from \cite{penzl00}. 
The system matrix $A\in  \mathbb{R}^{N\times N}$  arises from a spatial finite-difference discretization of the advection- diffusion equation
\begin{equation}\label{5.1}
\frac{\partial u}{\partial t} = \Delta u - 10x \frac{\partial u}{\partial x} - 100y \frac{\partial u}{\partial y}
\end{equation}
on the unit square domain $ \omega = (0, 1)^2$ with homogeneous Dirichlet boundary conditions. The matrices $B\in \mathbb{R}^{N\times 1}$ and $C\in \mathbb{R}^{1\times N}$ are defined as in \cite{Mena2018}.
The low-rank factor $L_0\in\mathbb{R}^{N\times 1} $ of the initial value $X(0)$ is generated randomly using MATLAB's \texttt{randn} function.

Each spatial dimension is discretized using $n_0 = 40$ uniformly spaced grid points, resulting in a system size of $N= n_0^2 = 1600$. The DRE is integrated over the time interval $[0, 0.1]$.

Fig.~\ref{fig4.1} reports the accuracy at the final time $t=0.1$ obtained by the four integrators as functions of both the number of time steps and the total CPU time. For the adaptive methods \texttt{exprb32} and \texttt{exprb43}, the plotted markers correspond to runs with prescribed tolerances $\texttt{tola}=\texttt{tolr}=10^{-1},10^{-2},10^{-3},10^{-4},10^{-5}$. The results for the non-adaptive methods \texttt{exprb2} and \texttt{exprb3} are obtained using fixed numbers of time steps $n=16,32,64,128,256$. The reference solution is computed by the fourth-order splitting scheme (\texttt{Split4}) with $n=10^{12}$ fixed time steps.
The results show that the adaptive integrators achieve higher accuracy with fewer time steps and reduced computational time compared with the non-adaptive counterparts. This confirms the effectiveness of the adaptive strategy in balancing accuracy and efficiency.

A main motivation for adaptive exponential integrators is their ability to resolve the transient phase more accurately. Fig.~\ref{fig4.2} (left) shows the evolution of the Frobenius norm of the numerical solutions over the interval $[0,0.1]$, together with the reference solution. The adaptive methods use $\texttt{tola}=\texttt{tolr}=10^{-5}$, whereas the non-adaptive methods employ a fixed number of time steps $n=256$. The right panel provides a magnified view of the transient regime. The results clearly indicate that the adaptive schemes compensate for the loss of accuracy observed in the non-adaptive methods during this phase by automatically selecting smaller time steps. 

Fig.~\ref{fig4.3} (left) displays the time-step evolution of the adaptive methods on a semi-logarithmic scale. The step size remains relatively small in the initial transient phase and then gradually increases as the solution approaches a steady-state regime. This behavior reflects the dynamical features of the problem and demonstrates the robustness of the adaptive controller. Fig.~\ref{fig4.3} (right) illustrates the accuracy achieved at the transient time $t=0.002$ as a function of the prescribed tolerance. The observed errors are in good agreement with the target tolerances, confirming the reliability of the error estimation strategy.
\end{example}

\begin{example}\label{exam2} In this example, we further assess the performance of the proposed adaptive exponential integrators \texttt{exprb32} and \texttt{exprb43}. Their performance is compared with the second-order Rosenbrock method \texttt{Ros2} and the fourth-order splitting method \texttt{Split4}. We consider two DRE generated from the benchmark datasets \texttt{flow} and \texttt{chip}, taken from \cite{morwiki18} and previously studied in \cite{Simoncini19}. These problems arise from real-world control applications and are commonly used to evaluate large-scale low-rank solvers.

The datasets provide the coefficient matrices $\widehat{E}, ~\widehat{A}\in \mathbb{R}^{N\times N}$,~ $\widehat{B}\in \mathbb{R}^{N\times q}$, and $\widehat{C}\in \mathbb{R}^{p\times N}$, originating from the linear control system
\begin{equation}\label{5.1}
\left\{
\begin{array}{l}
\widehat{E}~\widehat{x}~'(t)=\widehat{A}~\widehat{x}~(t)+\widehat{B}~u(t),\\
y(t)=\widehat{C}~\widehat{x}(t),
\end{array}
\right.
\end{equation}

The \texttt{flow} system has $N=9669$, $q=1$, and $p=5$, whereas the \texttt{chip} system has $N=20082$, $q=1$, and $p=5$. In both cases, the mass matrix $\widehat{E}$ is diagonal and positive definite. Introducing the change of variables $x(t)=\widehat{E}^{-1/2}\widehat{x}(t)$ leads to the standard form
\begin{equation}\label{5.2}
\left\{
\begin{array}{l}
x'(t)=Ax(t)+Bu(t),\\
y(t)=Cx(t),
\end{array}
\right.
\end{equation}
where
 \begin{equation*}
A=\widehat{E}^{-\frac{1}{2}}\widehat{A}~\widehat{E}^{-\frac{1}{2}}, \quad B=\widehat{E}^{-\frac{1}{2}}\widehat{B}, \quad C=\widehat{C}~\widehat{E}^{-\frac{1}{2}}.  
\end{equation*}

This transformation yields DRE of the form \eqref{1.1}  with low-rank structure \eqref{2b.1}.
Following \cite{Simoncini19}, the initial low-rank factors are chosen as $L_0=0$ for \texttt{flow} and $L_0=\sin(z)$ for \texttt{chip}, where $z\in\mathbb{R}^{N\times 1}$ has entries uniformly distributed in $[0,2\pi]$.

The two DRE are integrated over the time intervals $[0,0.001]$ and $[0,0.01]$, respectively. Tables~\ref{tab4.1} and~\ref{tab4.2} report the relative errors and CPU times obtained with \texttt{exprb32} and \texttt{exprb43} using $\texttt{tola}=\texttt{tolr}=10^{-5}$, and for \texttt{Ros2} and \texttt{Split4} using 256 fixed time steps. Again, the reference solutions are computed by \texttt{Split4} with $10^{12}$ time steps.

For the \texttt{flow} system, both adaptive methods achieve very high accuracy (down to $10^{-10}$) while requiring substantially less computational time than \texttt{Ros2}. Among the tested methods, \texttt{exprb32} provides the best overall efficiency. For the more challenging \texttt{chip} problem, the adaptive integrators maintain errors below the prescribed tolerance, although the observed accuracy is slightly reduced due to the strongly transient dynamics.

At the final time $t=0.01$, \texttt{exprb43} becomes less competitive in terms of CPU time because of its higher per-step cost, whereas \texttt{exprb32} consistently offers the best balance between accuracy and efficiency. Overall, the results demonstrate that adaptive exponential integrators are particularly effective for large-scale  DRE. They automatically adjust the time step to resolve rapid dynamics while reducing the computational effort during smoother phases. In contrast, splitting schemes remain efficient in short transient regimes but may require finer discretizations to maintain accuracy over longer time horizons.
\end{example}

\begin{table}[h]
\caption{The relative errors and the CPU times (in seconds) for the \texttt{flow} system of Example \ref{exam2}.}\label{tab4.1}
\begin{tabular*}{\textwidth}{@{\extracolsep\fill}lcccc}
\toprule%
 & \multicolumn{2}{@{}c@{}}{$t=0.001$}& \multicolumn{2}{@{}c@{}}{$t=0.01$} \\\cmidrule{2-3}\cmidrule{4-5} %
Methods &Relative errors &Times &Relative errors &Times  \\
\midrule
\texttt{Ros2}  &1.0780e-05&683.25   &8.0695e-06&1711.59 \\
\texttt{Split4}   &7.9061e-09&501.45  &3.5312e-06&1078.51   \\
\texttt{exprb32}    &3.0419e-10&104.72   &1.2822e-09&743.91  \\
\texttt{exprb43}    &3.0373e-10 &169.01 &1.2213e-09&1282.31 \\
\bottomrule
\end{tabular*} 
\end{table}

\begin{table}[h]
\caption{The relative errors and the CPU times (in seconds) for the \texttt{chip} system of Example \ref{exam2}.}\label{tab4.2}
\begin{tabular*}{\textwidth}{@{\extracolsep\fill}lcccc}
\toprule%
& \multicolumn{2}{@{}c@{}}{$t=0.001$}& \multicolumn{2}{@{}c@{}}{$t=0.01$} \\\cmidrule{2-3}\cmidrule{4-5}%
Methods &Relative errors &Times &Relative errors &Times  \\
\midrule
\texttt{Ros2}  &1.0278e-05&3793.76    &1.9188e-04&7732.85 \\
\texttt{Split4}   &2.0311e-08&443.78  &8.8084e-07&1266.35   \\
\texttt{exprb32}   &2.3180e-06&246.95   &1.9768e-06&942.40  \\
\texttt{exprb43}   &2.2696e-06&262.28 &1.1411e-06&1463.62 \\
\bottomrule
\end{tabular*} 
\end{table}

\section{Conclusion}\label{sec:5}
In this work, we have developed adaptive low-rank exponential Rosenbrock integrators for large-scale stiff differential Riccati equation. The proposed framework combines embedded error estimation with structure-preserving low-rank representations, allowing the time step to adjust automatically to the evolving solution dynamics while maintaining computational and memory efficiency. By performing all key operations in low-dimensional subspaces, the resulting methods achieve scalability with respect to the state dimension and effectively capture both transient and steady-state behaviors. Numerical experiments on standard benchmark problems demonstrate that the proposed adaptive integrators provide reliable error control and often improves both accuracy and computational efficiency compared with fixed-step low-rank integrators.


\bmhead{Funding}
 This work was partially supported by the Jilin Scientific and Technological Development Program (Grant No. YDZJ202501ZYTS635), the Natural Science Foundation of Jilin Province (Grant No. JJKH20241001KJ), and the National Natural Science Foundation of China (Grant No. 12371455).

\bmhead{Code and data Availability}
 The codes described in the manuscript are available at \url{https://github.com/lidping/adaptive-exprb-for-DRE.git}.

\section*{Declarations}
\bmhead{Ethics approval} Not applicable.
\bmhead{Competing interests} The authors declare to have no competing interests related to this work.

\bibliography{references}

\end{document}